\documentclass[aps,preprint,showkeys,amsmath,amssymb,
aps,
]{revtex4-1}

\usepackage{graphicx}
\usepackage{dcolumn}
\usepackage{bm}
\usepackage{color}
\usepackage[utf8]{inputenc}
\usepackage[toc,page]{appendix}

\usepackage{mathtools}
\usepackage{empheq}
\usepackage[skins,theorems]{tcolorbox} 
\tcbset{highlight math style={enhanced,
		colframe=black,colback=white,arc=4pt,boxrule=1pt}}

\usepackage{algpseudocode}
\usepackage{algorithm}
\usepackage[english]{babel}

\usepackage{etoolbox}
\usepackage{amsmath}

\newcommand{\algstrut}[1][\algruledefaultfactor]{\vrule width 0pt
	depth .25\baselineskip height #1\baselineskip\relax}
\newcommand*{\algrule}[1][\algorithmicindent]{\hspace*{0.5em}\vrule\algstrut
	\hspace*{\dimexpr#1+1.0em}}

\makeatletter
\newcount\ALG@printindent@tempcnta
\def\ALG@printindent{%
	\ifnum \theALG@nested>0
	\ifx\ALG@text\ALG@x@notext
	\else
	\unskip
	\ALG@printindent@tempcnta=1
	\loop
	\algrule[\csname ALG@ind@\the\ALG@printindent@tempcnta\endcsname]%
	\advance \ALG@printindent@tempcnta 1
	\ifnum \ALG@printindent@tempcnta<\numexpr\theALG@nested+1\relax
	\repeat
	\fi
	\fi
}%

\patchcmd{\ALG@doentity}{\noindent\hskip\ALG@tlm}{\ALG@printindent}{}{\errmessage{failed to patch}}

\AtBeginEnvironment{algorithmic}{\lineskip0pt}

\begin{document}

\date{\today}



\title{Beyond partial control: Controlling chaotic transients with the safety function}

\author{Rub{\'e}n Cape{\'a}ns}
\affiliation{Nonlinear Dynamics, Chaos and Complex Systems Group, Departamento de  F\'isica, Universidad Rey Juan Carlos, M\'ostoles, Madrid, Tulip\'an s/n, 28933, Spain}

\author{Miguel A.~F. Sanju{\'a}n}
\email{Corresponding author: miguel.sanjuan@urjc.es}
\affiliation{Nonlinear Dynamics, Chaos and Complex Systems Group, Departamento de  F\'isica, Universidad Rey Juan Carlos, M\'ostoles, Madrid, Tulip\'an s/n, 28933, Spain}

\begin{abstract}
	Partial control is a technique used in systems with transient chaos. The aim of this control method is to avoid the escape of the orbits from a region $Q$ of the phase space where the transient chaotic dynamics takes place.	This technique is based on finding a special subset of $Q$ called the safe set. The chaotic orbit can be sustained in the safe set with a minimum amount of control. In this work we develop a control strategy to gradually lead any chaotic orbit in $Q$ to the safe set by using the safety function. With the technique proposed here, the safe set can be converted into a global attractor of $Q$. 

\end{abstract}

\maketitle

\section{Introduction}

The control of chaotic systems developed in the late 20th century has been one of the main achievements in the field of nonlinear dynamics. The deterministic nature of these systems does not make them predictable. Due to the high sensitivity of the chaotic systems, the uncertainty in the forecast of the trajectories increases exponentially with time, making their control a challenging task. The OGY control method \cite {OGY} was the first to achieve the stabilization of an unstable periodic orbit embedded in a chaotic attractor, by applying small perturbations in a system parameter. Since then, a variety of control methods have been developed in different scenarios where chaotic dynamics is present.

One that has attracted interest in recent years is the chaotic transient behaviour. This situation usually arises  when due to the change of a parameter of the system, a chaotic attractor collides with his own boundary basin causing a boundary crisis,  giving the trajectories a path to escape towards an external attractor. Sooner or later almost all trajectories escape with the exception of a set called the chaotic saddle, where arbitrarily long transients can be found. This set is fractal and non-attracting and is the skeleton of the transient chaos behaviour \cite{TransientChaos2}.

In some situations, the chaotic behaviour is a desirable feature that is worth to preserve \cite{Thermal,Ecology}. However, if we have transient chaos, we need to apply some control strategy to sustain it. In this sense,  different control methods have been proposed to convert transient chaos into  permanent chaos.  Among the most important methods, we can refer to the methods proposed by Yang et al. \cite{Biological}, Schwartz and Triandaf \cite{Schwartz}, and Dhamala and Lai \cite{Dhamala}. With the same goal a new method called partial control has been proposed \cite{Asymptotic,Automatic}. The main difference of this method is that is conceived to deal with chaotic dynamics affected by noise. This actually makes a big difference,  since all real experiments are affected by noise. Even if this amount of noise is small, it should be taken into account since chaotic motion is an error amplifier, and small deviations at the beginning can ruin even the best control strategy. 

The intrinsic instability of the chaotic saddle, together with the action of noise creates a difficult scenario where keeping the control small might seem not realistic. However, the application of the partial control technique comes up with a surprising result. It was discovered that a region $Q$ with a chaotic saddle contains a special subset where chaotic trajectories can be sustained with a very small amount of control. This subset of $Q$ called the safe set, depends on the noise strength affecting the map and the available control intensity. The partial control technique is based on finding these special subsets, and has been successfully applied in many different problems \cite{Automatic, Ecology, Lorenz, Parametric}.The more surprising result is that the amount of control necessary to sustain the chaotic trajectory in the safe set is smaller than the amount of noise affecting it.

In all these previous works, partial control is limited to find a safe set in a certain region $Q$ of phase space and then sustain the chaotic trajectory inside the safe set. This approach implies that the initial condition of the orbit must already be in the safe set. However it may happen in some cases that the choice of the initial condition might be imposed by the problem.  For example, if we deal with the problem of putting a satellite into a particular orbit, the launch from Earth is an unavoidable initial condition that should be taken into account in the control design. For this reason, we will develop here a novel strategy to extend the control of the trajectories to any initial condition in $Q$, keeping the control as low as possible.

The paper is organized as follows. In Sect. 2, we briefly introduce the partial control method, and explain how to obtain the safe set from the safety function by using the slope-three tent map. In Sect. 3, we present an strategy going beyond the partial control method to gradually approach the safe set for any initial condition in $Q$. In Sect. 4, we analyze the variation of the safe set in function of the upper disturbance bound and the parameter $\mu$ of the tent map. Finally, we describe the main conclusions of our work.

\section{The partial control method}

Partial control is a control technique to sustain transient chaotic orbits affected by noise. The method is applied on a certain region $Q$ of phase space where orbits exhibit transient chaos that eventually escape from $Q$. It is assumed that the dynamics in $Q$ can be described with the following map:

\begin{equation}
\begin{array}{l}
q_{n+1}=f(q_n)+\xi_n+u_n, \hspace{1cm}\text{with} \hspace{0.4cm}
|\xi_n| \leq \xi_0, \hspace{0.3cm} |u_n| \leq u_0 < \xi_0,\\
\nonumber
\end{array}
\end{equation}
where $\xi$ is the noise (we call it disturbance) affecting the map and is considered to be bounded so that $|\xi_n|\leq\xi_0$.  The control term $u_n$ is also bounded so that $|u_n|\leq u_0$.

\begin{figure}
	\includegraphics [trim=0cm 0cm 0cm 0cm, clip=true,width=0.9\textwidth]{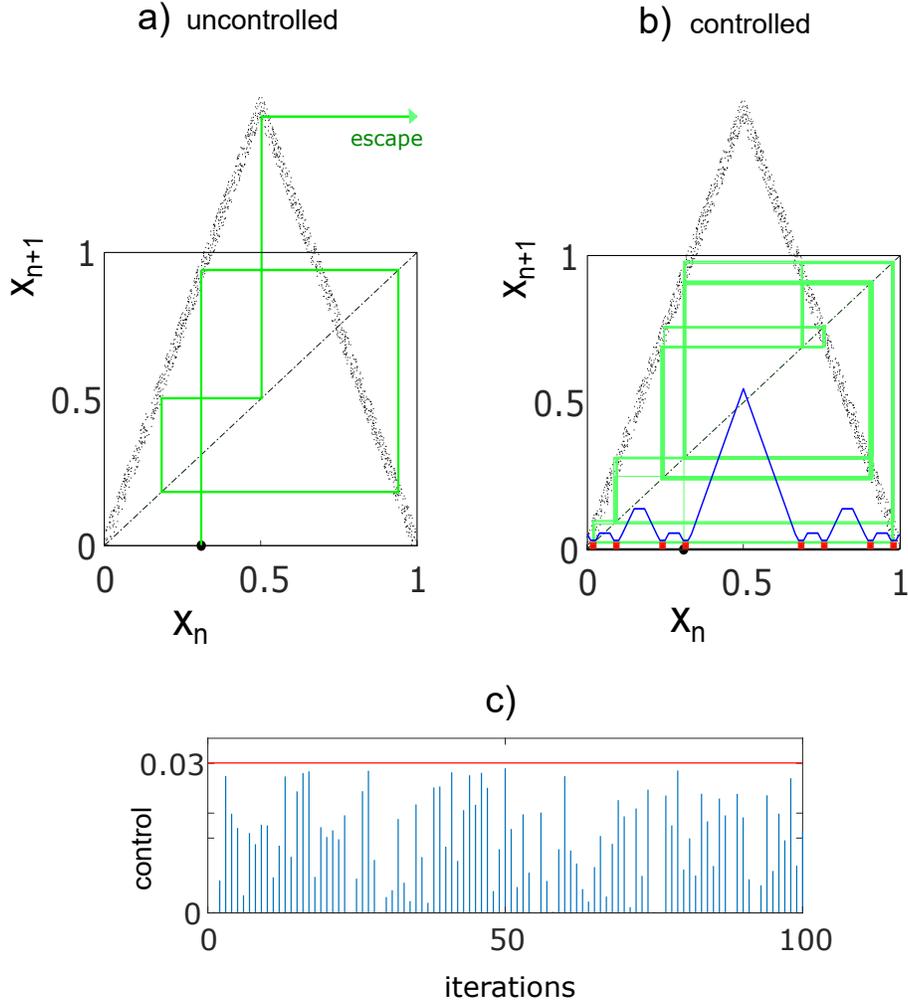}\\
	\centering
	\caption{\textbf{Partial control method.} The slope-three tent map is represented in this figure. The map is affected by a uniform disturbance bounded by $\xi_0=0.05$. The small dots help to visualize the intensity and distribution of the disturbance. (a) An uncontrolled orbit that escapes from the interval $Q=[-0,1]$ after a few iterations. (b) In blue, the safety function. It can be computed that the minimum value of the safety function is $0.03$. The minimal values of the safety function define the safe set, which is represented with the red pieces at the bottom. In green, a controlled orbit is shown. This controlled orbit starts in the point $x_0=0.3$ that belongs to the safe set. At every iteration of the map, the orbit is forced to pass through the safe set to remain forever in $Q=[0,1]$. (c) Controls $|u_n|\leq u_0$ applied during the first $100$ iterations of the map.}
	\label{1}
\end{figure}

The safe set is defined to be the subset of $Q$ where the controlled orbits can be sustained with $u_0<\xi_0$. This set can be obtained from a special function called the \textit{ safety function}, which has been developed in a previous work \cite{SafetyFun}. An algorithm to compute the safety function is described in an Appendix at the end. This function named $U$ is defined in each point $q\in Q$. The value $U(q)$ represents the minimum control bound that an orbit starting in $q$ needs to remain in $Q$ forever. For example, a value $U(q)=0.1$ means that a chaotic orbit with the initial condition $q$ can be sustained in $Q$ by applying each iteration a control $|u_n| \leq 0.1$.

The relation between the safe set and the safety function is the following. Given an upper disturbance bound affecting the map $\xi_0$, and an upper control bound $u_0<\xi_0$, the safe set corresponds to the points $q\in Q$ that satisfies  $U(q) \leq u_0$. The minimum  possible value $u_0$ is the minimum  value of the function $U$. Below this value, no safe set exist. In the following, we compute the safe sets with the minimum  value $u_0$.

To show a simple example of the partial control application, we will use here the well-known slope-three tent map. The equation of the map, including the disturbance term $\xi_n$ and the control term $u_n$, is given by: 

\begin{equation}
\label{}
x_{n+1} = \left\{
\begin{array}{ll}
\mu x_n + \xi_n +u_n      & \mathrm{\;for\ } x_n \le \frac{1}{2} \\
\mu(1-x_n) + \xi_n +u_n    & \mathrm{\;for\ } x_n > \frac{1}{2} \\
\end{array}
\right.
\end{equation}

\vspace{1cm}

This map with $\mu=3$, exhibits transient chaos in the interval $Q=[0,1]$ (see Fig.~\ref{1}a). In order to avoid the escape of the chaotic orbit, we compute first the corresponding safety function (see Fig.~\ref{1}b), where we consider an upper disturbance bound $\xi_0=0.05$. The used algorithm for its computation is detailed in~\cite{SafetyFun}. The safety function has $8$ minima with the value $u_0=0.03$.  The location of this minima defines the safe set, shown in Fig.~\ref{1}b by the small red pieces. 

Once the safe set is computed, the chaotic orbit can be controlled by forcing it to pass through the safe set at each iteration. Each control $|u_n|\leq 0.03$ is applied  with the knowledge of $f(q_n)+\xi_n$ so that $f(q_n)+\xi_n +u_n$ falls in the closest safe point. The control values corresponding to the first $100$ iterations is shown in Fig.~\ref{1}c.

\section{Beyond the partial control method}

As mentioned in the introduction, partial control applies only to initial points in the safe set, which is the subset of $Q$  where the upper control bound satisfies $u_0<\xi_0$. As shown before, orbits starting in this subset require the minimum control to remain in $Q$. 

\begin{figure}
	\includegraphics [trim=0cm -0.5cm 0cm 0cm, clip=true,width=0.9\textwidth]{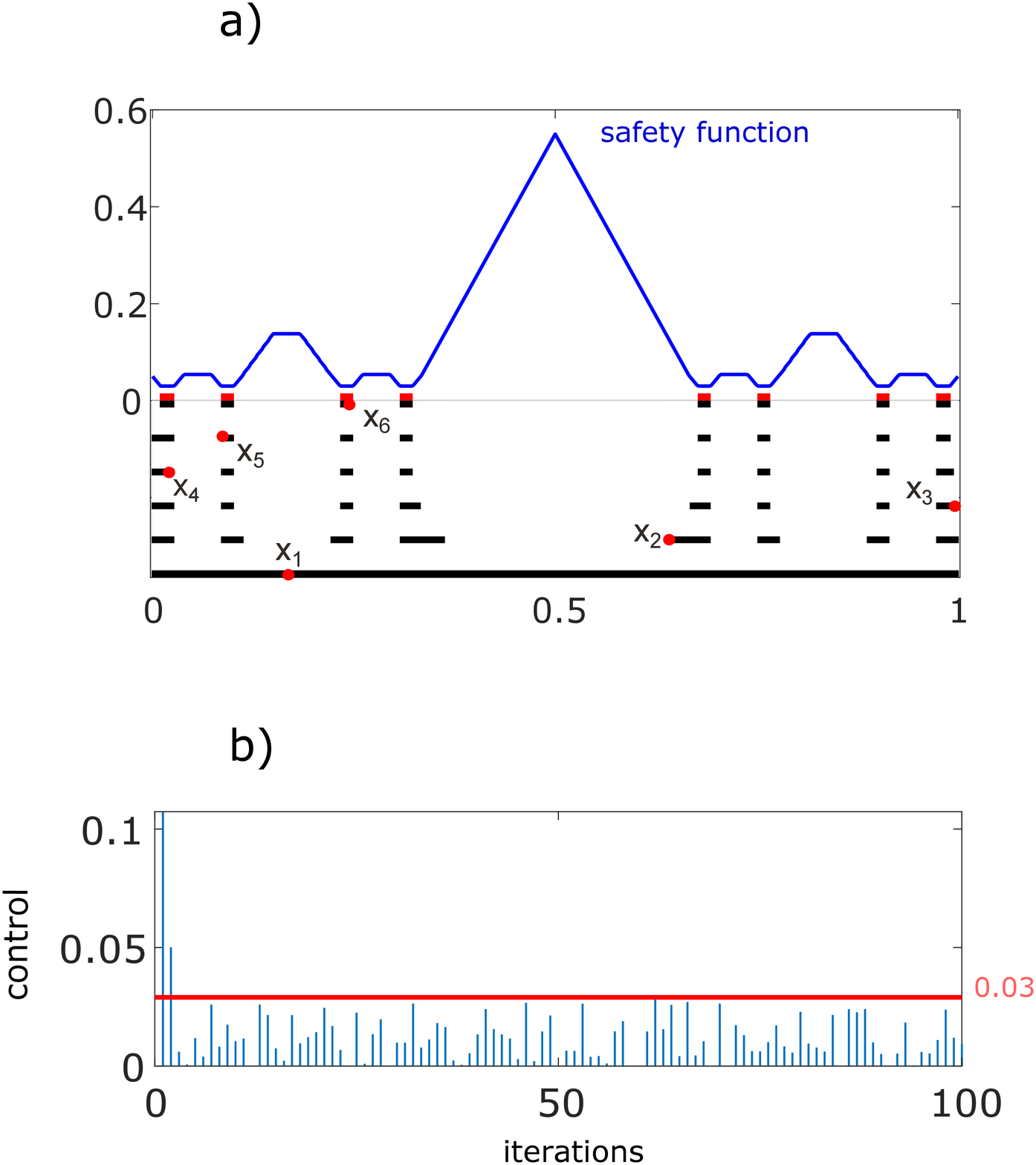}\\
	\centering
	\caption{\textbf{Controlling any initial condition in  $Q$}. (a) The black horizontal bar on the bottom represents a set of $1000$ initial conditions in $Q$. The remaining black pieces represent how these initial conditions gradually approach the safe set, shown in red. Any initial condition converges to the safe set in at most $6$ iterations. A particular controlled orbit is shown by the red dots. (b) The figure shows the first $100$ controls $|u_n|$ corresponding to the orbit shown in the figure, represented by the red dots. Note that the intensity of the controls decreases very quickly as the orbit approaches the safe set. Once the orbit enters into the safe set (in our example in at most $6$ iterations), the controls $|u_n|$ remain below the minimum control bound $0.03$, that is, the minimum value of the safety function.}
	\label{2}
\end{figure}

However, it may occurs that the initial condition of our orbit is outside the safe set. As the shape of the safety function shows (Fig.~\ref{2}a), initial conditions outside the safe set may require a very large upper control bound in comparison to the control needed in the safe set. The successive application of this large control can be prohibitive for the controller. Fortunately, we have observed that with a suitable strategy, the use of large controls is only needed at the very beginning. Broadly speaking, the  strategy consists of converting the safe set into a global  attractor of the orbits in the region $Q$. In this manner, the controlled orbits gradually approach the safe set, where the control needed is minimum.

Notice that, leading  the orbit to the safe set is equivalent to leading the orbit to the minima of the safety function $U$.  Therefore, given an initial condition $q_n\in Q$ with the corresponding value $U(q_n)$, the control $u_n$ is chosen so that the image $f(q_n)+\xi_n +u_n$ falls in a point $q_{n+1}$ with $U(q_{n+1}) \leq U(q_n)$. By repeating this strategy, the orbit gradually approaches to the minimum of the safety function, since $U(q_n) \geq U(q_{n+1}) \geq U(q_{n+2})..$. Eventually, the orbit enters into the safe set to remain forever in it.
 
However, this strategy still allows different choices to apply the control $u_n$ in each iteration, since the choice of the image $q_{n+1}$ with $U(q_{n+1}) \leq U(q_n)$ is not unique. Our proposal here is to select the control $u_n$ that minimizes the control bound of the orbit. This strategy can be implemented as follows:

 \begin{enumerate}
 	
 	\item Given an initial point $q_n\in Q$, evaluate the noisy image  $q^*=f(q_n)+\xi_n$.
 	\item Compute all the possible controls $|u_i|= |q_i-q^*|$ with $q_i\in Q$, being $i=1:N$ the  grid points in $Q$.
 	\item Among all the possible controls, apply the control $u_n = \min\limits_{1\leq i\leq N} \, \Big(\max \,\big(\,|u_i|, U(q_i)\,\big)\,\Big)$. The final point will be $q_{n+1}=f(q_n) + \xi_n + u_n =q^* + u_n$.
 	\item Repeat the algorithm with the new point $q_{n+1}$.
 	
 \end{enumerate}

\begin{figure}
	\includegraphics [trim=0cm -0.5cm 0cm 0cm, clip=true,width=0.9\textwidth]{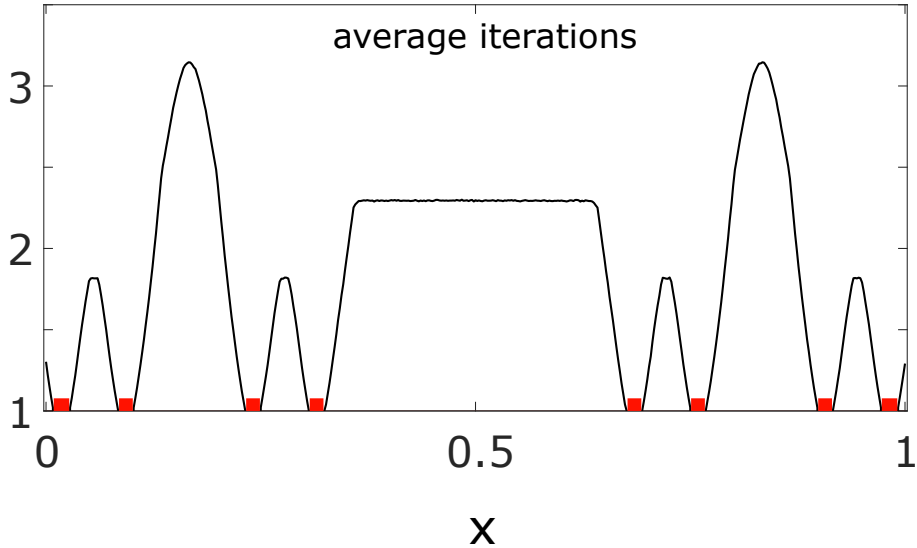}\\
	\centering
	\caption{\textbf{Average number of iterations to reach the safe set}. For each initial condition $q\in Q$, the average was done over $1000$ different orbits. The initial conditions farther from the safe set take more iterations to reach the safe set. The central region is an exception because this region maps directly outside the interval $Q$. See the explanation in the text.}
	\label{3}
\end{figure}

To show an example of this control strategy, we use again the slope-three tent map affected by a disturbance bound value $\xi_0=0.05$. The corresponding safety function is represented in Fig.~\ref{2}a. This function is the same as the one displayed in Fig.~\ref{1}b with a minimum value $u_0=0.03$. To show how any initial condition converges to the safe set,  we take a grid of $1000$ points in the interval $Q=[0,1]$. This set of initial conditions is represented with a black horizontal bar at the bottom of the Fig.~\ref{2}a. Then we apply the control strategy to the corresponding orbits. As shown, the initial set gradually converges to the safe set in at most $6$ iterations of the map. 

We also displayed in Fig.~\ref{2}a a particular controlled orbit, the orbit starting in $x_1=0.16$ (outside the safe set) and the next $5$ iterations after which the orbit reaches the safe set. The sequence of the first $100$ controls $|u_n|$ for this orbit is represented in Fig.~\ref{2}b. As shown, only the first controls applied are larger since the orbit starts outside the safe set. Once the orbit enters into the safe set, the controls remain below the minimum control bound $u_0=0.03$.

\begin{figure}
	\includegraphics [trim=0cm -0.5cm 0cm 0cm, clip=true,width=0.9\textwidth]{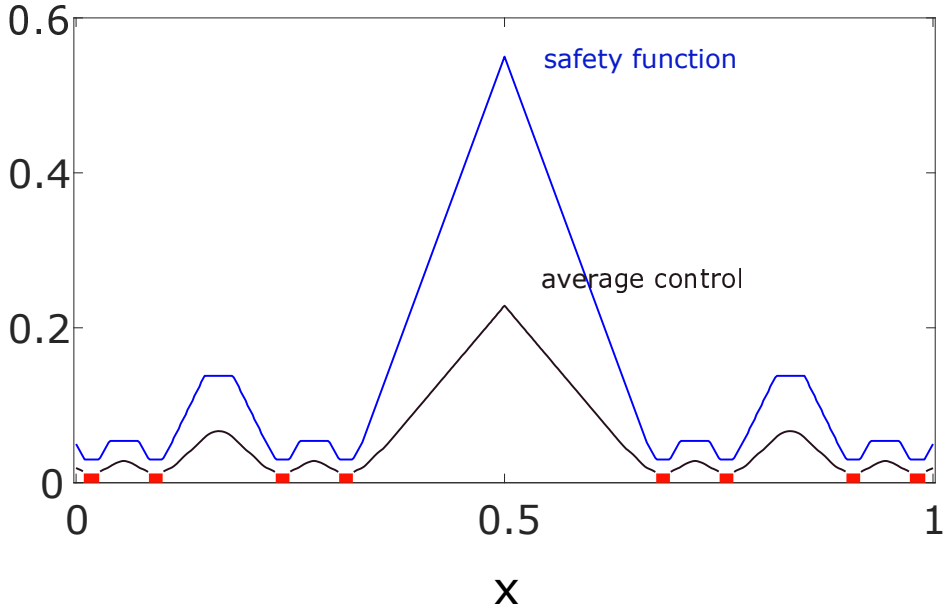}\\
	\centering
	\caption{\textbf{Average control per iteration to reach the safe set}. For each orbit, we compute the sum of the individual controls $|u_n|$ to reach the safe set divided by the number of iterations. For each initial point $q \in Q$, the average control has been computed for $1000$ different orbits. Note that the average control shape is similar to the safety function shape, showing that initial conditions that need a larger upper control bound also require a larger average control.}
	\label{4}
\end{figure}

Other quantity that could be interesting to learn about the rate of convergence towards the safe set, is the average number of iterations that an initial condition $q\in Q$ takes to reach the safe set, see Fig.~\ref{3}.  As expected, the initial conditions far from the safe set are the ones that need more iterations. The only discordant region is the central region (around $x=0.5$), where the average number of iterations is smaller than expected. The reason is because the central region maps directly outside the interval $Q=[0,1]$, in particular on the right side outside of $Q$.  Therefore, these points are directly reinserted into the right side inside of $Q$ (applying a large control), where the safety function takes low values. As a consequence, the initial conditions placed in the center of $Q$,  quickly reach the safe set. The counterpart is the large control used. 

To measure the size of the control needed, we represent in Fig.~\ref{4}, the average control that each initial condition $q\in Q$ needs to reach the safe set. For the computation of this average control, we only count the number of iterations that the orbit takes to reach the safe set. As it can be seen, the shape of the average control resembles the shape of the safety function. This shows that the initial conditions that need a larger upper control bound to remain in $Q$, also require a larger average control to converge to the safe set. The farther an initial condition is from the safe set, the more control it will need. 

\section{Variation of the safe set with $\xi_0$ and $\mu$}

In the previous section we have seen how the safe set can be transformed in a global attractor of the initial conditions in $Q$. Once the trajectories enter into the safe set, the upper control bound $u_0$ required is minimum. 

As shown, the safe set is the subset of $Q$ where the safety function is minimum. This function strongly depends on the upper disturbance bound $\xi_0$ affecting the map and the parameter value of the system. As a consequence, the location of the safe set and the corresponding control $u_0$ may change.

\begin{figure}
	\includegraphics [trim=0cm 0cm 0cm 0cm, clip=true,width=0.9\textwidth]{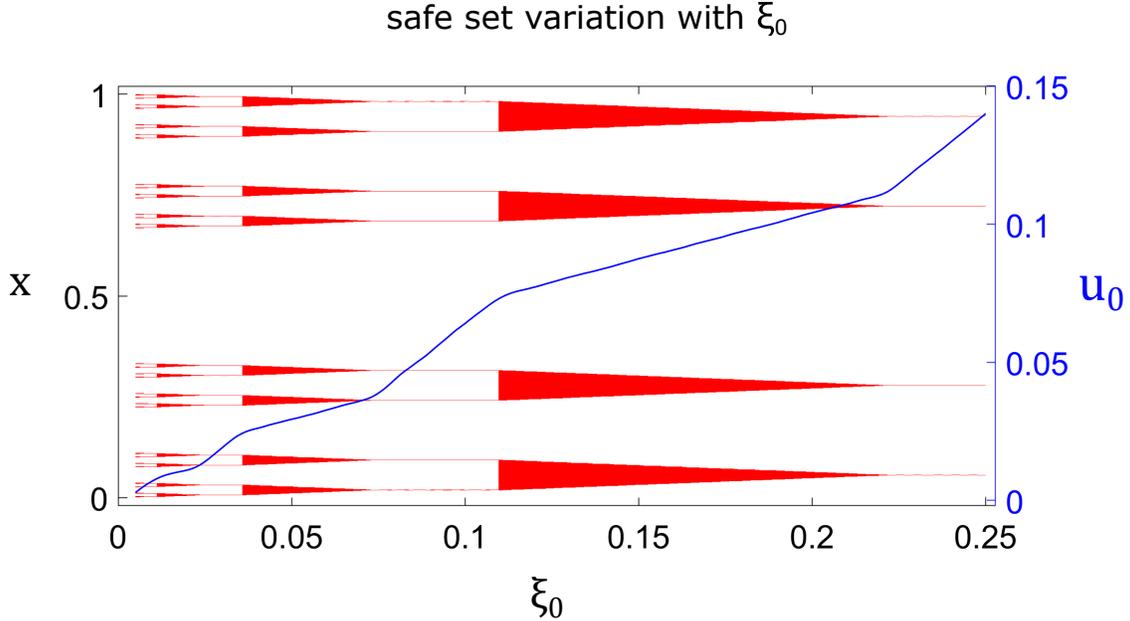}\\
	\centering
	\caption{\textbf{Safe set variation with $\xi_0$.} In this plot we represent in red different safe sets computed for different values $\xi_0=[0.005,0,25]$. The parameter $\mu=3$ of the tent map is kept fixed. In blue, the upper control bound $u_0$ corresponding to each safe set. Notice that the relation $u_0$/$\xi_0$ is almost constant and smaller than $1$, as a matter of fact about 0.6.}
	\label{5}
\end{figure}

To visualize this variation, we display in Fig.~\ref{5} the different safe sets in red computed for the slope-three tent map, that is with the parameter $\mu=3$, and different disturbance values in the range $\xi_0=[0.005,0.25]$. We represent here the interval $Q=[0,1]$ on the left y-axis. For a given value $\xi_0$, the safe set corresponds to the intersection of the vertical axis located at this value with the red figure. On the other hand, we plot on the right y-axis, the $u_0$ control value in blue associated to each safe set at a given value $\xi_0$. As it can be seen, from left to right, the safe sets structure resembles the Cantor structure characteristic of the chaotic saddle, the topological object responsible for the transient chaotic dynamics. Lower values of $\xi_0$ result in safe sets with many small pieces, whereas larger $\xi_0$ values result in safe sets with only a few pieces. Note that each piece of the safe set is a region of $Q$ where the safety function has a minimum.  It is also remarkable how the thickness of these pieces suddenly changes (see for example the change in $\xi_0=0.11$).  With respect to  the variation of $u_0$, we can see that larger $\xi_0$ disturbance values require larger $u_0$ control values as expected. It is noteworthy to observe that the relation $u_0$/$\xi_0$ is almost constant and always smaller than $1$, as a matter of fact about $0.6$.

\begin{figure}
	\includegraphics [trim=0cm 0cm 0cm 0cm, clip=true,width=0.9\textwidth]{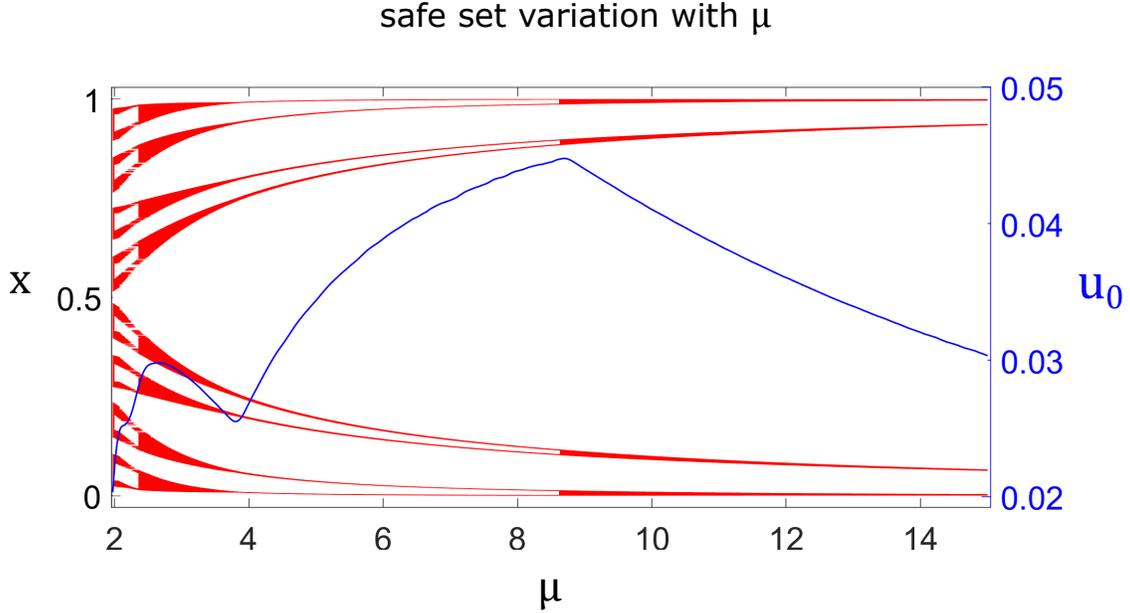}\\
	\centering
	\caption{\textbf{Safe set variation with $\mu$.} In this plot, different safe sets in red have been computed for different values in the range $\mu=[2,15]$. The upper disturbance bound $\xi0=0.05$ is kept fixed. In blue, the upper control bound $u_0$ corresponding to each safe set. The value $u_0$ is lower where the separation $d$ between nearby pieces of the safe set is approximately the same as the disturbance value $\xi_0=0.05$.}
	\label{6}
\end{figure}

Now, we fix the upper disturbance bound $xi_0=0.05$ and we have computed different safe sets in the parameter range $\mu=[2,15]$ of the tent map. This is shown in Fig.~\ref{6}.  As in Fig.~\ref{5}, the left y-axis is the interval $Q=[0,1]$ and the right y-axis represents the $u_0$ control value corresponding to each safe set. On the left side of the figure we can observe that for $\mu=2$, the control $u_0$ is not zero. This is because, the action of the disturbance in the tent map, allows the orbits to escape from $Q$ even for parameter values below $\mu=2$. We can also observe that, for increasing values of $\mu$, the number of pieces of the safe set decrease, with sudden changes in $\mu=2.35$ and $\mu=8.67$. 

However, what is more surprising in Fig.~\ref{6} is how the upper control bound $u_0$ varies. This control increases and decreases without a clear pattern showing that, larger parameter values of $\mu$ do not necessary implies larger controls $u_0$. What we have observed is that the average separation $d$ between nearby pieces of the safe set, has an important influence in the value $u_0$ required. The value $u_0$ tends to be lower in the safe sets for which $d \approx \xi_0=0.05$. For example, from left to right in the figure, for $\mu=2.6$ the separation $d$ between nearby pieces is about $0.1$, which is far from $\xi_0$, and the control $u_0$ has a local maximum. For  $\mu=3.8$ however, the separation $d$ between nearby pieces is about $0.05$, and the control $u_0$ has a local minimum. Then, the control increases again until it reaches a maximum for $\mu=3.8$. Here, the separation $d$ between nearby pieces of the safe set is approximately $0.015$, also far from $\xi_0$. At this point, the $8$ pieces of the safe set collapses in $4$ pieces and the separation $d$ changes to approximately $0.1$, again far from $\xi_0$.  Finally, for parameter values $\mu>3.8$ the separation $d$ decreases, taking values closer and closer to $\xi_0$. The $u_0$ value in this region also decreases until the last safe set plotted ($\mu=15$), where the separation $d$ is about $0.06$.

The reason for this relation between $d$ and $\xi_0$ can be understood if we take into account that, roughly speaking, each piece of the safe set maps in the middle between other nearby pieces. The implications of this behavior can be seen with the next simple example. Imagine we have a map affected by $\xi_0=1$ and the safe set only has two pieces, made of two points. The left piece, for example, is placed in $q_l=-0.5$ and the right piece in $q_r=0.5$. Notice that $q_r-q_l=\xi_0$ and the middle point between $q_r$ and $q_l$ is $0$. Therefore, $f(q_r)=f(q_l)=0$. After adding the disturbance and the control, we have that $q_{n+1}=f(q_n)+\xi_n+u_n =0+\xi_n+u_n$ with $|\xi_n| \leq \xi_0=1$ and $|u_n| \leq u_0$. Now, we are interested to find the minimum upper control bound $u_0$ required to control the orbit, that is, to put the point $q_{n+1}$ again in the safe set, that is the set containing the points $q_r$ or $q_l$. To find  $u_0$, we analyze the worst case of disturbance $\xi_n$ which are $\xi_n=0$ or $|\xi_n|=1$. Thus, the control required for these cases is $0.5$ and therefore the upper control bound will be $u_0=0.5$.

Imagine now that varying a parameter of the map, the locations of $q_r$ and $q_l$ change so that $q_r-q_l\neq\xi_0$. It is clear that if $q_r-q_l < \xi_0$ the worst disturbance is  $|\xi_n|=1$ and if  $q_r-q_l > \xi_0$  the worst disturbance is $\xi_n=0$. In both cases, the upper control bound required will be $u_0>0.5$. Therefore the separation $d$ that minimizes the value $u_0$ is $d=\xi_0=1$. This approach can be extended to safe sets with more than two pieces as shown in Fig.~\ref{6}. This shows that $u_0$ tends to be lower in the safe sets for which the separation between nearby pieces satisfy $d \approx\xi_0=0.05$. This is merely a heuristic explanation lacking rigor, as it goes without saying that features such as the shape of the map can also influence the determination of $u_0$.

\section{Conclusions}

The partial control technique is able to sustain a transient chaotic orbit by means of the computation of the safe set. This special set is the subset of $Q$ where the safety function has a minimum. Orbits with initial conditions inside the safe set need a minimum control to remain in $Q$.  However, orbits with initial conditions outside the safe set may need a larger control to remain in $Q$.  Here, we present a different strategy to greatly reduce these large controls. This control strategy consists of gradually approaching the orbit to the safe set. As a result, the safe set is an attractor of any initial condition in $Q$.

We have illustrated the application of this method in the slope-three tent map. For this map, we have shown how the initial conditions in the interval $Q=[0,1]$ converges to the safe set. Both, the average iterations to converge and the average control applied, strongly depend on how far the initial condition is from the safe set. In any case, after a few iterations, the orbit enters into the safe set where it can be sustained with minimum control $u_0$.

We have also analyzed how the safe set and the corresponding control $u_0$ change depending on the disturbance $\xi_0$ and the parameter $\mu$ of the tent map. Our results reveal that the larger the value $\xi_0$, the larger the value $u_0$. In contrast, the variation of $u_0$ with $\mu$ is more complex, showing that the separation $d$ between nearby pieces of the safe set plays an important role.

Finally, we want to emphasize that, although the chaotic map used here was one-dimensional for simplicity, this control technique can be easily extended to higher dimensional maps. The only requirement is the extra computation effort to compute safety functions in higher dimensions.


\begin{algorithm}[H]
	\renewcommand{\thealgorithm}{}
	\floatname{algorithm}{}
	\caption{\textbf{Appendix: The safety function algorithm}}\label{euclid}
	
	\begin{algorithmic}[0]
	\Statex
	\Statex

	Map with the index notation:\\ $~q_{n+1}=f(q_n) + \xi_n +u_n$ $~~~\rightarrow~~~$  $q[j]=\,f\big(q[i]\big) +\xi[s]+ u\,[i,s,j]$.\\
	$i\equiv$ index of the starting point $q[i],~ i=1:N$ \,\,\,\,\,\, with $N=$ total number of grid points in $Q$.\\
	$s\equiv$ index of the disturbance $\xi[s],~ s=1:M$  \,\,\,\,\,\,\,\, with  $M=$ number of disturbed images.\\
	$j\equiv$ index of the arrival point $q[j],~ j=1:N$.\\\\
	
	\textbf{Computation of the safety function $U$}:

	\State- Initially set $\;U_0\,[j]=0, \; \forall j=1:N, \; \; \; \; k=0.$
	
	\hspace{4cm}\While {$U_{k+1}\neq U_k$}
	
	\For{$i=1$ to  $N\;$}
	\For{$s=1$ to  $M\;$}
	\For{$j=1$ to  $N\;$}
	\State $\;u\,[i,s,j]= \Big|\,f\big(q[i],\xi[s]\big)- \;q\,[j]\,\Big|$
	\State
	\State $\;u^*\,[i,s,j]=\max\limits_{j}\,\big(\,u[i,s,j],\, U_k[j]\,\big)$
	\EndFor
	\State $\;u^{**}\,[i,s]=\min\limits_{j}\,\big(\,u^*\,[i,s,j]\,\big)$
	\EndFor
	\State $U_{k+1}[i]=\max\limits_{s} \big(\,u^{**}\,[i,s]\,\big)$
	\EndFor
	\State  $k=k+1$\normalsize
	\EndWhile
	\Statex

		\State {{Compact form of the while loop: $U_{k+1}[i]=\max\limits_{1\leq s\leq M}\Big(\min\limits_{1\leq j\leq N}\big(\max\limits_{j}\, (\,u[i,s,j],\, U_k[j]\,)\,\big)\,\Big)$}} \\
		\State {Note: The $u\,[i,s,j]$ values remain unchanged every iteration of the while loop so only it is necessary to compute them once and save them.}

	\end{algorithmic}
\end{algorithm}


\begin{acknowledgments}
This work was supported by the Spanish State Research Agency (AEI) and the European Regional Development Fund (FEDER) under Project No. PID2019-105554GB-100.
\end{acknowledgments}

\section*{Compliance with ethical standards}
\textbf{Conflict of interest} The authors declare that they have no conflict of interest concerning the publication of this manuscript.

\section*{Data availability}
All data generated or analysed during this study have been obtained through numerical simulations that are included in this article

\end{document}